\documentclass[12pt, reqno]{amsart}
\usepackage{amsfonts}
\usepackage{graphicx}
\usepackage{amsmath, amscd, amssymb}
\usepackage[latin5]{inputenc}

\begin{document}

\begin{center}
\thispagestyle{empty} 
\markboth{\textbf{ }}
{\textbf{Complex variable approach to analysis of a FDE}}

\noindent {\Large \textbf{Complex variable approach to analysis of a fractional differential equation in the real line \\[0pt]}}

\bigskip

\textbf{M\"ufit \c{S}an }
\end{center}

\vspace{.03cm}

\begin{center}
\noindent \textit{Department of Mathematics, Faculty of Science, Çank{\i}r{\i} Karatekin University, Tr-18100, Çank{\i}r{\i}, Turkey\\[0pt]
e-mail: mufitsan@karatekin.edu.tr } \\[0pt]
\end{center}

\medskip

\vspace{.4cm}

\noindent \textbf{Abstract.} The first aim of this work is to establish a Peano's type existence theorem for an initial value problem involving complex fractional derivative, and the second is, as a consequence of this theorem, to give a partial answer to the local existence of the continuous solution for the following problem 
\begin{equation*}
\begin{cases} 
&D_{x}^{q}u(x) = f\big(x,u(x)\big),      \\  
&u(0)=b, \ \ \ (b\neq 0).   \\
\end{cases} 
\end{equation*}
Moreover, in the special cases of considered problem, we investigate the some geometric properties of the solutions.  \\ 


\medskip


\bigskip

\section{\noindent Introduction and Motivation}

In the theory of ordinary differential equation, Peano existence theorem is practical and important because one can easily check the local existence of the solution to a differential equation. Such type existence theorem in the theory of fractional differential equation was first given by Laksmikantham and Vatsala in \cite{Laks} by using Tonelli's approach for the following problem 
\begin{equation*}
\begin{cases} 
&D_{x}^{q}u(x) = f\big(x,u(x)\big),    \tag{1.1}  \\  
&u(0)=b,\ \       \\
\end{cases} 
\end{equation*}
where $D_{x}^{q}$ is well-known Riemann-Liouville fractional derivative in the real line, $b$ is a real number and $f\in C([0,T]\times\mathbb{R},\mathbb{R}).$ However, as indicated in  \cite{Zhang},  the problem (1.1) with $b\neq 0$ is not a suitable problem. For this reason, this problem was considered with a different initial data or with a modified equation by some researchers, for example \cite{Bale2}, \cite{Zhang}.  Nevertheless, if the initial condition in (1.1) is taken as homogenous, one can see that the problem (1.1) has meaning provided that $f\in C([0,T]\times\mathbb{R},\mathbb{R}).$ In any case, the problem (1.1) with the non-homogenous initial condition ($b\neq 0$) can be suitable, if the nonlinear function $f$ in this problem satisfies the following conditions: \\ 

(i) Let the functions $f(x,y)$ and $x^{q}f(x,y)$ ($0<q<1$) be continuous on $(0,T]\times\mathbb{R}$ and $[0,T]\times\mathbb{R},$ respectively. \\

(ii) $x^{q}f(x,b)\big|_{x=0}=b/\Gamma(1-q).$ \\

It is known from the study of Delbosco and Rodino \cite{Del} that if the condition (i) holds, then the equation (1.1) is equivalent to the following integral equation: 
$$
u\left(x\right) \ = \frac{1}{\Gamma \left( q \right) }\  \int_{0}^{x}\frac{\ f\left( \zeta ,u\left( \zeta \right) \right) }{\left(x-\zeta
\right) ^{1-q }}d\zeta, \eqno{(1.2)}
$$
where $\Gamma$ is the well-known Gamma function.\\

Now, if the initial condition with $b\neq 0$ in (1.1) is considered, then the condition (ii) is necessary for the problem to be suitable. Otherwise, the following contradiction arises:   
$$b=u(0)=\lim_{x\rightarrow 0^+}u(x)=\frac{1}{\Gamma(q)}\lim_{x\rightarrow 0^+} \int_{0}^{1}\frac{(xt)^{q}\ f\left(xt,u\left(xt\right) \right) }{t^{q}\left( 1-t\right)^{1-q }}dt\neq b.$$   

It needs to note that the existence and uniqueness of continuous solution to this problem without initial condition was proved in Theorem 3.5 in \cite{Del}  by posing an Lipschitz-type condition in addition to (i). However, the local existence of continuous solution to this problem under only the condition (i) is still an open problem. \\

In this work, we give a partial answer to this open problem by establishing a Peano-type existence theorem for the complex version of the problem (1.1): 
\begin{equation*}
\begin{cases} 
&D_{z}^{q}u(z) = f\big(z,u(z)\big),  \ \ z\in\mathbb{U}_{R},  \tag{1.3}  \\  
&u(0)=b,    \\
\end{cases} 
\end{equation*}
where $\mathbb{U}_{R}:\big\{z\in\mathbb{C}:\left|z\right|<R\big\}$ is a open disc in the complex plane $\mathbb{C},$  $b\in\mathbb{C}$ and $D^{q}$ is the complex fractional derivative operator which is a direct extension of the well-known Riemann-Liouville derivative in the real line and given in \cite{Samko},\cite{San}. \\

The existence of solution to the problem (1.3) is investigated in the Banach space $\mathcal{B}(\mathbb{U}_{R}):=C(\overline{\mathbb{U}}_{R})\cap \mathcal{A},$ where $\mathbb{U}_{R}:=\big\{z\in\mathbb{C}:\left|z\right|<R\big\},$ $C(\overline{\mathbb{U}}_{R})$ is denoted the space of continuous functions on $\overline{\mathbb{U}}_{R},$ and  $\mathcal{A}$ is the space of analytic functions on $\mathbb{U}_{R},$ provided that the nonlinear function $f$ satisfies the following two conditions:  \\

(iii) $f(z,t)$ is a multivalued function on $\overline{\mathbb{U}}_{R}\times \mathbb{C}$  such that $z^{q}f(z,t)$ is
analytic on $\mathbb{U}_{R} \times \mathbb{C}$ and continuous on $\overline{\mathbb{U}}_{R}\times \mathbb{C}.$ 
\\ 
  
(iv) $z^{q}f(z,b)\big|_{z=0}=b/\Gamma(1-q).$   \\ 

As a result of this existence theorem (See Corollary 2.6 and Example 2.7), it is possible to obtain the existence of continuous solution of the problem (1.1) in the following way:  If the nonlinear function $f(x,t)$ satisfying the conditions (i)-(ii)  can be extended to the the function $f(z,t)$ satisfying the condition (iii)-(iv), then as a consequence of this existence theorem we can say that there exists a solution $u\in\mathcal{B}(\mathbb{U}_{R})$ of (1.3) such that real part of which is a continuous solution of the problem (1.1). From this point of view, we think that this existence theorem may contribute to the area of fractional differential equation in the real line. Moreover, in the proof of this theorem we use a new technique related to Schwarz Lemma for analytic function of one variable and two variables varying on generally different closed balls in $\mathbb{C},$ which we prove in the next section. From this aspect, this existence result is a nice consequence of dealing with analytic functions. \\

On the other hand, as shown by Diethelm in \cite{Diethelm}, the existence of real analytic solution of the problem 1.1 with Caputo-derivative $D_{x}^{q}$ is a rare event when $f$ is analytic. It can be said that the same result is valid for the problem 1.1 and 1.3  when $f$ is analytic. \\

Furthermore,  the complex fractional integral and derivative operators have been widely used in the complex analysis and have many applications in the univalent function theory (see, for example, \cite{Chen},\cite{Irmak},\cite{Sri}). We here show that the solutions of this problem with the special cases of $f$ are univalent or starlike on $\mathbb{U}$ by using the results of Noshiro-Warschawski in \cite{Goodman} and Mocanu in \cite{Mocanu}. Such investigations were previously made for the solutions to the ordinary differential equations in \cite{Owa2},\cite{Saitoh},\cite{San2}.

\section {Preleminaries and Main Result}

In this section, we intend to present a Peano's type existence theorem for the problem (1.3) involving fractional derivative $D_{z}^{q},$ and therefore the fractional integral $I_{z}^{q},$ which are defined as follows (see \cite{Samko}):     \\

\noindent \textbf{Definition 2.1.} Let the function $u(z)$ defined on a certain domain of complex plane containing the points $0$ and $z.$ Then, the fractional integral and derivative of order $q$ $\big(0<q <1\big)$ of $u(z)$ are defined, respectively, by 
\begin{equation}
I_{z}^{q }u\left( z\right) :=\frac{1}{\Gamma \left( q \right) }%
\ \int_{0}^{z}\frac{u\left( \zeta \right) }{\left( z-\zeta \right)
^{1-q }}d\zeta,  \tag{2.1}
\end{equation}
and 
\begin{equation}
D_{z}^{q }u\left( z\right) =\frac{1}{\Gamma \left( 1-q \right) 
}\frac{d}{dz}\int_{0}^{z}\frac{u\left( \zeta \right) }{\left( z-\zeta
\right) ^{q }}d\zeta   \tag{2.2}
\end{equation}
where the integrations are along the straight line interval connecting points $0$ and $z$ as a rule, and with the principal value
$$(z-\zeta)^{{1-q}}=\left|z-\zeta\right|^{1-q}e^{i(1-q)\arg z}, \  \ \ \arg z\in (-\pi,\pi]. $$ 

The definitions mentioned above are the direct generalization of the well-known Riemann-Liouville fractional integral and derivative in the real line. Since the integration in these definitions are over the line segment connecting points $0$ and $z,$ we can say that 
$$D_{x}^{q }u(x)= D_{z}^{q }(\Re\big(u\left( z\right)\big)) \ \text{on} \ (0,R) $$
for any $R>0$ and for $u(x),$ real part of $u\left( z\right).$ \\

Under the certain conditions, the compositional relations for these operations  
$$D_{z}^{q }I_{z}^{q }u(z)=u(z) \ \  \text{and} \ \ I_{z}^{q }D_{z}^{q}u(z)=u(z)$$
were given and proved in Lemma 3.2 and Remark 3.3 in \cite{San}.\\

By help of these relations, the following equivalence was revealed there to the investigate the existence of the solution in  $\mathcal{B}^{b}(\mathbb{U}_{R}):=\left\{u\in\mathcal{B}(\mathbb{U}_{R}): u(0)=b\right\}$ of the problem (1.3). \\ 

\noindent\textbf{Lemma 2.1.} Under the conditions (iii)-(iv), the problem (1.3) is equivalent to the following Volterra-type integral equation 
\begin{equation*}
u\left(z\right)=\frac{1}{\Gamma \left(q\right) }\int_{0}^{z}\frac{f\left(\zeta,u\left(\zeta\right)\right)}{\left( z-\zeta
\right) ^{1-q}}d\zeta, \ \ z\in\overline{\mathbb{U}}_{R}, \eqno{(2.3)}
\end{equation*}
that is every solution of (1.3) is also a solution of (2.3) and vice versa.\\ 

The existence theorems in the ordinary and fractional differential equations (for example, \cite{Del},\cite{Yu}, \cite{Zeidler}) were proved by means of Schauder's fixed point theorem. However, making only use of this theorem is not enough to prove the local existence of desired solution of the integral equation (2.3) for any $f$ satisfying the conditions (iii) and (iv). Let us see this in the following: \\
 
\noindent\textbf{Remark 2.2.} For simplicity, let us take $b=0$ in the problem (1.3) with $f$ satisfying the conditions (iii) and (iv). Then, for the fixed positive real numbers $R,$ $r$ there exists a $M$ such that 
$$\sup_{(z,t)\in\overline{\mathbb{U}}_{R} \times \overline{\mathbb{U}}_{r}}\left|z^{q}f(z,t)\right|\leq M, \ \  (M>0),$$    
where $\mathbb{U}_{R}:=\big\{z\in\mathbb{C}: \left|z\right|<r \big\}$ and $\mathbb{U}_{r}:=\big\{t\in\mathbb{C}: \left|t\right|<r \big\}$ are the open discs. Then, we have for the operator $T$ defined as (2.6)
\begin{align*}
\left|Tu\left(z\right)\right|\leq\frac{1}{\Gamma\left(q\right)}\int_{0}^{\left|z\right|}\frac{\left|\zeta^{q}f(\zeta,u(\zeta))\right|}{\left|\zeta\right|^{q} \left|z-\zeta \right|^{1-q}}d\zeta\leq M\Gamma(1-q)\leq r \ \ \text{on} \ \ \overline{\mathbb{U}}_{R_0}
\end{align*}
for any $R_0\leq R,$ if  $M\Gamma(1-q)\leq r.$ That is, we can mention about the existence of solution in $u\in\mathcal{B}^{0}(\mathbb{U}_{R}),$ only if $f$ satisfies $M\Gamma(1-q)\leq r$ in addition to the conditions (iii) and (iv). We meet this same result  when we investigate the local continuous solution of the equation in (1.1) (or the problem (1.1)) under the condition (i) (or the conditions (i) and (ii)). For this reason, this problem is open. \\

As a consequence of the explanation in Remark 2.2, in the proof of existence theorem we also use a new technique related to Schwarz Lemma for analytic function of one variable (See  \cite{Dettman}, \cite{Pon}, \cite{San}) and Schwarz Lemma for analytic function of two variables, which is given as follows: \\

\noindent\textbf{Lemma 2.3.}  Let $R,$ $r$ be fixed positive real numbers and $b$ be a fixed complex number. Let also $\mathbb{U}^{b}_{r}:=\big\{t\in\mathbb{C}: \left|t-b\right|\leq r \big\}$ and 
$g:\overline{\mathbb{U}}_{R} \times \overline{\mathbb{U}}^{b}_{r} \longrightarrow \mathbb{C}$ 
be analytic on $ \mathbb{U}_{R}\times\mathbb{U}^{b}_{r},$ continuous on $\overline{\mathbb{U}}_{R} \times 
\overline{\mathbb{U}}^{b}_{r}$ with $g(0,b)=0$ and bounded with $\left|g(z,t)\right|\leq M$ ($M>0$) for all 
$(z,t)\in\overline{\mathbb{U}}_{R} \times \overline{\mathbb{U}}^{b}_{r}.$ 
Then,
\begin{align*}
\left|g(z,t)\right|\leq M\max\left(\frac{\left|z\right|}{R},\frac{\left|t-b\right|}{r} \right) \tag{2.4}
\end{align*}
for all $(z,t)\in \overline{\mathbb{U}}_{R} \times \overline{\mathbb{U}}^{b}_{r}.$\\

\textit{Proof.} Let $\xi=(\xi_{1},\xi_{2})$ be fixed complex number on $\partial(\mathbb{U}_{R})\times\partial(\mathbb{U}^{b}_{r})$
and we define a function $\Phi_{\xi}$  by the following form:
\begin{align*}
\Phi_{\xi}:\overline{\mathbb{U}}_{R}\rightarrow \overline{\mathbb{U}}_{R} \times \overline{\mathbb{U}}^{b}_{r} \ \ , \ \ \Phi_{\xi}(\eta):=\left(\frac{\xi_{1}}{R}\eta,\frac{(\xi_{2}-b)}{R}\eta+b\right).
\end{align*}
Clearly, from the the hypotheses of this theorem, $g\circ\Phi_{\xi}$ is then analytic on $ \mathbb{U}_{R},$
continuous on $\overline{\mathbb{U}}_{R}$ and also bounded with $\left|(g\circ\Phi_{\xi})(\eta)\right|\leq M $ 
for all $\eta\in \overline{\mathbb{U}}_{R}.$  

Now, let a function $\psi(\eta)$ be defined by 
$$\psi(\eta):=\frac{(g\circ\Phi_{\xi})(\eta)}{\eta}, \ \ \ \big( \eta\in \overline{\mathbb{U}}_{R}\big).$$  
Because of $(g\circ\Phi_{\xi})(0)=0,$  the function $\psi$ has a removable singularity at the point $\eta=0.$  Moreover, $\psi$ is 
continuous both on $\partial(\mathbb{U}_{R})$ and at $\eta=0,$ since $\psi(0)=(g\circ \Phi_{\xi})'(0)$ and also 
$\psi$ is analytic function on $\mathbb{U}_{R}-\left\{0\right\}.$  It follows from Cauchy integral formula,  
$\left|\psi(0)\right|=\left|(g\circ \Phi_{\xi})'(0)\right|\leq \frac{M}{R}$ is then obtained and from the hypothesis of this theorem, 
$\left|\psi(z)\right| \leq \frac{M}{R}$ on  $\partial(\mathbb{U}_{R})$ is also obtained. In view of these two results and by applying  
maximum modulus theorem to the function $\psi,$ we easily get that $\left|(g\circ \Phi_{\xi})(\eta)\right|\leq \frac{M\left|\eta\right|}{R}$  on $\overline{\mathbb{U}}_{R},$ which immediately yields that inequality (2.4).   

Furthermore, the equality (2.4) holds for all $(z,t)\in \overline{\mathbb{U}}_{R}\times \overline{\mathbb{U}}^{b}_{r},$ since, for any 
$(z,t)\in \overline{\mathbb{U}}_{R}\times \overline{\mathbb{U}}_{r},$ there exists a fixed number  
$\xi^{*}:=\left(\xi^{*}_{1},\xi^{*}_{2}\right)\in\partial(\mathbb{U}_{R})\times\partial(\mathbb{U}^{b}_{r})$ and also
a function $\Phi_\xi^{*}:\overline{\mathbb{U}}_{R}\rightarrow \overline{\mathbb{U}}_{R} \times \overline{\mathbb{U}}^{b}_{r}$ 
such that 
\begin{align*}
(z,t)=\left(\frac{\xi^{*}_{1}}{R}\eta,\frac{\xi^{*}_{2}}{R}\eta\right)=\Phi_{\xi^{*}}(\eta).
\end{align*} 
This completes the proof. \\

Now, we are ready to give the Peano's type existence theorem which is our first aim:  \\

\noindent\textbf{Theorem 2.4.}  Let $q$ be fixed in $(0,1)$ and let the hypotheses (iii)-(iv) be satisfied. Furthermore, for the fixed positive real numbers $R,$ $r$ let $$\sup_{(z,t)\in\overline{\mathbb{U}}_{R} \times \overline{\mathbb{U}}^{b}_{r}}\left|z^{q}f(z,t)-\frac{b}{\Gamma(1-q)}\right|\leq M, \ \  (M>0).$$ 
Then, the initial value problem (1.3) possesses at least one solution $u$ in $\mathcal{B}^{b}(\mathbb{U}_{R_{0}}),$ where 
\begin{align*}
R_{0}:=
\begin{cases}
\ \ \ \   R&,  \ \text{if} \ \ M\Gamma(2-q)\leq r,\\ 
\smallskip
\frac{rR}{M\Gamma(2-q)}&, \ \text{if} \  \ r< M\Gamma(2-q).\\  
\end{cases}
\  \   \  \  \  \ \tag{2.5}
\end{align*}
               
\textit{Proof.}  As indicated in Lemma 2.1, the problem (1.3) is equivalent to the integral equation (2.3). Let us now define a operator $T$ as follows: 
$$Tu\left(z\right)=\frac{1}{\Gamma \left( q \right) }\int_{0}^{z}\frac{f\left( \zeta ,u\left(\zeta\right)\right)}{\left( z-\zeta
\right)^{1-q}}d\zeta, \ \ z\in\overline{\mathbb{U}}_{R}. \eqno{(2.6)} $$
It is clear that $T:\mathcal{B}^{b}(\mathbb{U}_{R_0}) \ \to \ \mathcal{B}^{b}(\mathbb{U}_{R_0}).$ Hence, the fixed points of the operator $T$ in $\mathcal{B}^{b}(\mathbb{U}_{R_0})$ coincide the solutions of the problem (1.3). Thus, our aim is to prove the existence of
the fixed points of the operator $T$ by means of Schauder's fixed point theorem. For the proof, we first show that the following inclusion  
\begin{align*}
T(B_{r})\subseteq B_{r}, \tag{2.7} 
\end{align*}
where $B_{r}$ is the appropriate closed, convex and bounded subset  of $\mathcal{B}^{b}(\mathbb{U}_{R_0}).$ 
The choice of this subset differs greatly depending on the value of $z^{q}f(z,t).$ This choice is made in the following way: let us take $g(z,t):=z^{q}f(z,t)-\frac{b}{\Gamma(1-q)}.$ Since the hypotheses of Lemma 2.3 are satisfied, we have
$$\left|z^{q}f(z,t)-\frac{b}{\Gamma(1-q)}\right|\leq M\max\left(\frac{\left|z\right|}{R},\frac{\left|t-b\right|}{r} \right) $$ 
for all $(z,t)\in \overline{\mathbb{U}}_{R} \times \overline{\mathbb{U}}^{b}_{r}.$ That is, there exist two case as follows: \\
(a) $\frac{\left|z\right|}{R}:= \max\left(\frac{\left|z\right|}{R},\frac{\left|t-b\right|}{r} \right) $ or, \\
(b) $\frac{\left|t-b\right|}{r}:= \max\left(\frac{\left|z\right|}{R},\frac{\left|t-b\right|}{r} \right).$ \\
For both cases above, we define a $B_{r}$ such that (2.7) is satisfied. When case (a) holds, let   
$$B_{r}= \big\{u\in\mathcal{B}^{b}(\mathbb{U}_{R_0}):\sup_{z\in\overline{\mathbb{U}}_{R_0}} \left|u(z)-b\right|\leq r\big\} $$ 
with $r<M\Gamma(2-q).$ \\

For any $u\in B_{r},$ from (2.6) we have
\begin{align*}
\left|Tu\left(z\right)-b\right|&=\left|\frac{1}{\Gamma \left(q\right) }\int_{0}^{z}\frac{f\left(\zeta,u\left(\zeta\right)\right)}{\left( z-\zeta
\right) ^{1-q}}d\zeta-b\right|\\
&\leq\frac{1}{\Gamma\left(q\right)}\int_{0}^{\left|z\right|}\frac{\left|f\left( \zeta ,u\left(\zeta\right)\right)-\zeta^{-q}\frac{b}{\Gamma(1-q)}\right|}{\left|z-\zeta\right|^{1-q}}\left|d\zeta\right|,   \\ 
&=\frac{1}{\Gamma\left(q\right)}\int_{0}^{\left|z\right|}\frac{\left|\zeta^{q}f\left( \zeta ,u\left(\zeta\right)\right)-\frac{b}{\Gamma(1-q)}\right|}{\left|\zeta\right|^{q}\left|z-\zeta\right|^{1-q}}\left|d\zeta\right|,  \\
&\leq\frac{M}{\Gamma\left(q\right)}\int_{0}^{\left|z\right|}\frac{\frac{\left|\zeta\right|}{R}}{\left|\zeta\right|^{q}\left|z-\zeta\right|^{1-q}}\left|d\zeta\right|,
\end{align*}
since it is supposed that case (a) holds. By using change of variable $\zeta=zt$ in the last inequality we obtain  
\begin{align*}
\left|Tu\left(z\right)-b\right|\leq\sup_{z\in\overline{\mathbb{U}}_{R_0}}\frac{M\Gamma\left(2-q\right)\left|z\right|}{R}\leq r,
\end{align*}
because of the definition of $R_0$ in (2.5). Hence, the inclusion in (2.7) is fulfilled. \\

Now, we suppose that case (b) holds and let us define 
$$B_{r}= \big\{u\in\mathcal{B}^{b}(\mathbb{U}_{R}):\sup_{z\in\overline{\mathbb{U}}_{R}}\left|u(z)-b\right|\leq r\big\} $$ 
with $M\Gamma(2-q)\leq r.$
By using same way in the above equations we have 
\begin{align*}
\left|Tu\left(z\right)-b\right|\leq\frac{M}{\Gamma\left(q\right)}\int_{0}^{\left|z\right|}\frac{\frac{\left|u(\zeta)-b\right|}{r}}{\left|\zeta\right|^{q}\left|z-\zeta\right|^{1-q}}\left|d\zeta\right| \ \text{for all} \ z\in\overline{\mathbb{U}}_{R}. \tag{2.8}
\end{align*}
On the other hand, from Schwarz Lemma of one variable we can write 
$$\frac{\left|u(z)-b\right|}{r}\leq \frac{r\left|z\right|}{rR}=\frac{\left|z\right|}{R} \ \ \text{ for all} \ \ z\in\overline{\mathbb{U}}_{R}.$$
If we use this inequality in (2.8) and the fact $M\Gamma(2-q)\leq r,$ then we obtain
\begin{align*}
\left|Tu\left(z\right)-b\right|\leq\frac{M}{\Gamma\left(q\right)}\int_{0}^{\left|z\right|}\frac{\frac{\left|\zeta\right|}{R}}{\left|\zeta\right|^{q}\left|z-\zeta\right|^{1-q}}\left|d\zeta\right|\leq M\Gamma(2-q)\leq r,  
\end{align*}
which shows that (2.7) holds for the case (b). \\

Let us now show that $T$ is a continuous operator on $B_{r}$ (for both cases). It is supposed that $\left \{u_{n}\right \}_{n=1}^{\infty}\subset B_{r}$ is a sequence with $u_{n}\rightarrow u$ in $\mathcal{B}^{b}(\mathbb{U}_{R_0})$ as $n\rightarrow \infty.$ Then, it follows that $u_{n}$ converges uniformly to $u\in B_{r}.$ By using this and the uniform continuity of $z^{q}f(z,t)$ on $\overline{\mathbb{U}}_{R}\times \overline{\mathbb{U}}^{b}_{r},$ we conclude that
\begin{align*}
&\left|Tu_{n}(z) -Tu(z)\right|=\left \vert \frac{1}{\Gamma \left(
q \right) }\  \int_{0}^{z}\frac{ \left[f\left( \zeta ,u_{n}(\zeta
)\right) -f\left( \zeta ,u(\zeta )\right)\right] }{\left( z-\zeta \right)
^{1-q}}d\zeta \right \vert \\
&\ \ \leq \frac{1}{\Gamma \left(q \right)} \int_{0}^{1}\frac{ \left|(\xi z)^{q}f\left(\xi z,u_{n}(\xi z)\right) -(\xi z)^{q}f\left( \xi z ,u(\xi z)\right)\right|}{\xi^{q}\left( 1-\xi \right) ^{1-q}} d\xi \rightarrow 0
\end{align*}%
as $n\rightarrow \infty.$ \\

Now let us prove  that $T(B_{r})$ is an equicontinuous set of $\mathcal{B}^{b}(\mathbb{U}_{R_0}).$ 
For any $u\in B_{r},$  $z^{q}f(z,u(z))$ is uniformly continuous on $\overline{\mathbb{U}}_{R}.$ Therefore, for given $\epsilon>0$ there exists a $%
\delta=\delta(\epsilon)>0$ such that 
$$\left \vert z_{1}^{q}f(z_{1},u(z_{1}))-z_{2}^{q}f(z_{2},u(z_{2}))\right \vert <\frac{\epsilon}{\Gamma(1-q)},$$
for all $z_{1},z_{2}\in \overline{\mathbb{U}}_{R}$ satisfying $\left|z_{1}-z_{2}\right|<\delta.$ From here, one can conclude that 
\begin{align*}
&\big|Tu\left(z_{1}\right)-Tu\left(z_{2}\right) \big| \\
&\ \ \ \leq \frac{1}{\Gamma \left(q \right)}\int_{0}^{1} \frac{\left \vert
(\xi z_{1})^{q}f\left(\xi z_{1},u\left( \xi z_{1}\right) \right) -(\xi
z_{2})^{q}f\left( \xi z_{2},u\left( \xi z_{2}\right) \right) \right
\vert}{\xi^{q} \left( 1-\xi \right) ^{1-q }} d\xi \  \\
&\ \ \ <\Gamma(1-q)\frac{\epsilon}{\Gamma(1-q)}=\epsilon,
\end{align*}
since $\left|\xi z_{1}-\xi z_{2}\right|<\delta.$ So, the desired is obtained.  \\

Therefore, as a consequence of Schauder fixed point theorem, it can be said that the operator $T$ has at least one fixed point in $\mathcal{B}^{b}(\mathbb{U}_{R_0})$ for a $R_0$ given in (2.5), which is also a solution of the problem (1.3). This completes the proof. \\

Before we give a corollary of Theorem 2.5 dealing with the existence of the solution for the initial value problem (1.1) with $b\neq 0,$ let us make the following preparations: We suppose that the nonlinear function $f(x,t)$ in (1.1) satisfies the conditions (i) and (ii). Let us extend the function $f(x,y)$  to the function $f(z,t)$ by writing $z$ and $t$ instead of $x$ and $y,$ respectively, for $(z,t)\in\overline{\mathbb{U}}_{R}\times \mathbb{C}.$ Therefore, we obtain a new function $f(z,t).$  Now, we are ready to give the following corollary:  \\

\noindent\textbf{Corollary 2.6.}  If this new function $f(z,t)$ mentioned above fulfills the conditions (iii)-(iv), and 
$$f(x,y)=\Re\big(f(z,t)\big)$$ 
holds, then the problem (1.1) with $b\neq 0$ has at least one continuous solution on the interval $[0,R_0]$ ($R_0$ is as in (2.5)), which is real part of the solution of the problem (1.3) under the conditions (iii)-(iv).   \\

\textit{Proof.}  If the conditions of this corollary are considered, one can say that, as a result of the previous theorem, there exists a solution $u(z)$ in $\mathcal{B}^{b}(\mathbb{U}_{R_0})$ for a $R_0$ given in (2.5) of the following  problem
\begin{equation*}
\begin{cases} 
&D_{z}^{q}u(z) = f\big(z,u(z)\big),      \\  
&u(0)=b, \ \ (b\in\mathbb{R}-\{0\}).    \\
\end{cases}
\end{equation*}
Thus, one can easily deduce from above that $\Re\big(u(z)\big)$ satisfies the above problem on $z\in[0,R].$  \\

Moreover, from the assumptions, one can write: \\
$$x^{q}f(x,y)=z^{q}\Re\big(f(z,t)\big) \ \text{for all} \ z\in[0,R].$$ \\
Therefore, it is obtained that $\Re\big(u(z)\big)$ satisfies the problem:
\begin{equation*}
\begin{cases} 
&D_{x}^{q}\left(\Re\big(u(z)\big)\right)= f\big(x,\Re\big(u(z)\big)\big),      \\ 
&u(0)=b, \ \ (b\in\mathbb{R}-\{0\}),  \\
\end{cases}
\end{equation*}
that is, $\Re\big(u(z)\big)$ is a continuous solution of the problem (1.1). \\

\noindent\textbf{Example 2.7.}  If we take $f(x,y)=\frac{x^{-q}}{\Gamma(1-q)}\left(y+\frac{q}{1-q}x\right)$ in the problem (1.1), then it is clear that $f(z,t)=\frac{z^{-q}}{\Gamma(1-q)}\left(t+\frac{q}{1-q}z\right)$ fulfills the conditions (iii)-(iv). As a result of Theorem 2.5, we can say that there exist a solution in $\mathcal{B}^{b}(\mathbb{U}_{R_0})$ of the problem (1.3) and this solution is $u(z)=b+z.$ Moreover, real part of $u(z),$ $u(x)=b+x$ satisfies the problem (1.1). \\    

The uniqueness of desired solution of the problem (1.3) was previously proved Theorem 3.8 in \cite{San}. In the following, we give two results related to the geometric properties of the unique solution on the unit disc $\mathbb{U}$ for the problem (1.3) with $u(0)=0$ and the function $f$ satisfying certain conditions, by using the Noshiro-Warschawski theorem in \cite{Goodman} and a result of Mocanu in \cite{Mocanu}. \\

\noindent\textbf{Theorem 2.8.}  Let $f(z,t):=z^{-q}h(z)$ with $h\in\mathcal{B}^{0}(\mathbb{U})$ in the problem (1.3). Then the problem (1.3) with $u(0)=0$ has a unique univalent solution, if $h$ satisfies the following conditions: \\

(i) $h'(0)=\frac{1}{\Gamma(2-q)},$ \\

(ii) $\Re\big(e^{i\beta}h'(z)\big)>0$ is satisfied for some real $\beta$ and for all $z\in\mathbb{U}.$\\
 
\textit{Proof.} By virtue of Theorem 3.8 in \cite{San}, one can say that there exists a unique solution $u$ in $\mathcal{B}^{0}(\mathbb{U})$ of this problem and, this solution is given by 
$$u\left(z\right)=\frac{1}{\Gamma \left(q\right)}\int_{0}^{z}\frac{\zeta^{-q}h(\zeta)}{\left( z-\zeta
\right) ^{1-q }}d\zeta. \eqno{(2.9)}$$
From here, it follows, by applying the change of variable $\zeta=zt,$ and by differentiating the both sides, that 
\begin{equation*}
u'\left(z\right)=\frac{1}{\Gamma \left(q\right)}\int_{0}^{1}t^{1-q}(1-t)^{q-1}h'(zt)dt, \ \ (\forall z\in\mathbb{U}).  
\end{equation*} 

Using the condition (ii) in the above equality, it is easily seen that  $\Re\big(e^{i\beta}u'(z)\big)>0$ for some real $\beta$ and for all $z\in\mathbb{U}.$  Moreover, if we use the condition (i) in the equality (2.9) by taking $z\rightarrow 0,$ then we have $u'(0)=1.$ Therefore, from the Noshiro-Warschawski theorem, the solution $u$ of the considered problem is univalent on $\mathbb{U}.$ 
\\    

\noindent\textbf{Theorem 2.9.}  Let $f(z,t):=z^{-q}h(z)$ with $h\in\mathcal{B}^{0}(\mathbb{U})$ in the problem (1.3). Then this problem has a unique starlike solution, if $h$ satisfies the condition (i) in Theorem 2.8 with the following one: \\

(i*) There exists $M$ with $0<M\leq \frac{\sqrt{20}}{5}$ such that the following inequality is fulfilled: 
\begin{align*}
\sup_{z\in\overline{\mathbb{U}}} \left|\Gamma(1-q)h(z)-\frac{z}{1-q}\right|\leq M.
\end{align*}  
 
\textit{Proof.}  As a result of Theorem 3.8 in \cite{San}, it can be said that there exists a unique solution $u\in\mathcal{B}^{0}(\mathbb{U})$ of the considered problem and it is given by (2.9). From (2.9), one can write for all $z\in\overline{\mathbb{U}}:$
\begin{align*}
\left|u(z)-z\right|&=\frac{1}{\Gamma(q)}\left|\int_{0}^{z}\frac{\left[\zeta^{-q}h(\zeta)-\frac{\zeta^{1-q}}{\Gamma(2-q)}\right]}{(z-\zeta)^{1-q}}d\zeta \right|  \\
&\leq\sup_{z\in\overline{\mathbb{U}}} \left|\Gamma(1-q)h(z)-\frac{z}{1-q}\right|. 
\end{align*}
By using the condition (i*) in the last inequality, it is easily obtained that $\sup_{z\in\overline{\mathbb{U}}}\left|u(z)-z\right|\leq M$ where $0<M\leq \frac{\sqrt{20}}{5}.$ From the Cauchy integral formula, we have 
$$\left|u'(z)-1\right|=\frac{1}{2\pi}\left|\int_{\left|\zeta\right|=1}\frac{u(\zeta)-\zeta}{(\zeta-z)^{2}}d\zeta\right|\leq M$$
for all $z\in\overline{\mathbb{U}}.$ 

Therefore, if the result of Mocanu in \cite{Mocanu} is considered in the last inequality, then it is seen that the considered solution $u$ is starlike on 
$\mathbb{U}.$

\end{document}